\documentclass[12pt]{article}
\usepackage{amssymb}

\def\appendix{\par}  
\def\cc{{\mathfrak c}}
\def\tt{{\mathfrak t}}

\def\proof{\par\noindent Proof\par\noindent}

\def\qed{\par\noindent QED\par}
\def\rmor{\mbox{ or }}

\def\rmand{\mbox{ and }}
\def\om{\omega}
\def\al{\alpha}

\def\del{\delta}
\def\diag{\Delta}
\def\infsets{{[\om]^\om}}
\def\sq{\subseteq}
\def\sm{\setminus}
\def\sig{\sigma}
\def\uu{{\mathcal U}}
\def\vv{{\mathcal V}}
\def\ww{{\mathcal W}}
\def\bb{{\mathcal B}}

\def\CC{{\mathcal C}}

\newtheorem{theorem}{Theorem}
\newtheorem{lemma}[theorem]{Lemma}
\newtheorem{question}[theorem]{Question}
\newtheorem{prop}[theorem]{Proposition}

\begin{document}

\begin{center}
{On squares of spaces and $F_\sig$-sets}
\end{center}

\begin{flushright}
Arnold W. Miller\footnote{
Thanks to the Fields Institute for Research in Mathematical Sciences
at the University of Toronto for their support during the time this paper
was written and to Juris Steprans who directed the special program in set
theory and analysis.  Also I want to thank Gary Gruenhage for telling
me about the problem and his translation
of the topology into properties of sets of reals.
}
\end{flushright}

\begin{quote}
Abstract: We show that the continuum hypothesis implies there exists
a Lindel\"of space $X$ such that $X^2$ is the union of two metrizable
subspaces but $X$ is not metrizable.  This gives a consistent solution
to a problem of Balogh, Gruenhage, and Tkachuk. The main lemma is
that assuming the continuum hypothesis there exist disjoint sets
of reals $X$ and $Y$ such that $X$ is Borel concentrated on $Y$,
i.e., for any Borel set $B$ if $Y\sq B$ then $X\sm B$ is countable,
but $X^2\sm\diag$ is relatively $F_\sig$ in $X^2\cup Y^2$.
\end{quote}

In Balogh, Gruenhage, and Tkachuk \cite{bgt} the following question
is asked:

\bigskip
Question 4.1.  Let $X$ be a  
regular paracompact space $X$ such that $X\times X$ is the union of two
metrizable subspaces.  Must $X$ be metrizable?
What if $X$ is Lindel\"of?

\bigskip
\begin{theorem}\label{one}
Assume the continuum hypothesis. Then
there exists a nonmetrizable regular Lindel\"of space 
$X$ such that $X^2$ is the union of 
two metrizable subspaces.
\end{theorem}

We first prove the following Lemma.  

\begin{lemma}
(CH) There are
uncountable disjoint sets $X,Y\sq 2^\om$
such that 
\begin{enumerate}
 
 \item $X$ is Borel concentrated on
 $Y$, i.e. every Borel set in $2^\om$ containing $Y$ contains all but countably many
 elements of $X$,

 \item $Y^2\sm\diag$ is $F_\sig$ in $X^2\cup Y^2$, and
 
 \item $X^2\sm\diag$ is $F_\sig$ in $X^2\cup Y^2$.

\end{enumerate}
Here $\diag=\{(x,x): x\in 2^\om\}$.
\end{lemma}

\proof
We identify the Cantor space $2^\om$ with the
power set of $\om$, $P(\om)$.  We use
$\infsets$ to stand for the infinite subsets of
$\om$.   Define for $y\in \infsets$
$$[y]^{*\om}=\{x\in[\om]^\om: x\sq^* y\}$$
where $\sq^*$ stands
for inclusion mod finite.
Let $\langle B_\al:\al<\om_1\rangle$ be all Borel subsets
of $[\om]^\om$.  We construct $y_\al$ for $\al<\om_1$
so that
\begin{enumerate}
\item $\al<\beta$ implies $y_\beta\sq^* y_\al$ and
$y_\beta\not=^* y_\al$ and
\item either $y_{\al}\notin B_\al$ or
$[y_{\al}]^{*\om}\sq B_\al$.
\end{enumerate}

These conditions are easy to get. Given $y_\beta$ for $\beta<\al$
and $B_\al$ let $y\in \infsets$ be arbitrary with 
$y\sq^* y_\beta$ but $y_\beta\not=^* y$ 
for each $\beta<\al$. 
If  $[y]^{*\om}$ is a subset of $B_\al$,
then simply take $y_{\al}\in [y]^{*\om}\sm B_\al$,
otherwise take $y_{\al}=y$.

Let
$$X=\{y_\al\sm y_{\al+1}:\al<\om_1\} \rmand  Y=\{y_\al:\al<\om_1\}$$

Iff $B$ is any Borel set containing $Y$, then 
choose $\al$ so that $B=B_\al$.  At stage $\al$ of
the construction it must have been that
$[y_{\al}]^{*\om}\sq B_\al$.  But this
means that $x_\beta\in B_\al$ for all $\beta\geq\al$. So $X$ is
Borel concentrated on $Y$. 

If we define 
$$F=\{(u,v)\in P(\om)\times P(\om): (u\sq^*v \rmor v\sq^*u)\rmand u\not=v\}$$
Then $F$ is an $F_\sig$ set and
 $$F\cap (X^2\cup Y^2) = (Y^2\sm\diag)$$
Also if we define 
 $$H=\{(u,v)\in P(\om)\times P(\om): u\cap v=^*\emptyset\}$$
Then $H$ is an $F_\sig$ set and
$$H\cap (X^2\cup Y^2) = (X^2\sm\diag)$$

This finishes the proof of the Lemma. 
\qed

Now define the following Michael-line like topology.  Suppose
that $M$ is a topological space and $X\sq M$.  Then $M(X)$
is the topological space on the same set but with the following
topology.  For $x\in X$
we make $x$ an isolated point, i.e., add $\{x\}$ to the topology
of $M(X)$.  For any point $y\in M\sm X$ neighborhoods
in $M$ form a neighborhood basis for $y$ in $M(X)$.
It is easy to see that $M(X)$ is regular for any regular
space $M$ and subset $X\sq M$.

The following is exercise 5.5.2 from Engelking \cite{engel}:

\begin{prop}\label{prop}
Suppose $M$ is a metric space and $X\sq M$.  Then $M(X)$ is metrizable
iff $X$ is an $F_\sig$ set in $M$.
\end{prop}

Our example is $M(X)$ where $X$ and $Y$ are from the
Lemma and $M=X\cup Y$ has its usual
(separable metric) topology as a subspace of $2^\om$.  
It follows from the Proposition that $M(X)$ is not metrizable. 

\bigskip\noindent
$M(X)$ is a Lindel\"of space.  Take any open cover $\uu$ of $M(X)$. 
Open sets in $M(X)$ have the form  
$U \cup Z$ where $U$ is open in
$M$ and $Z\sq X$ is arbitrary.
Then since $Y$ has its
standard topology, countably many elements of $\uu$ will cover $Y$,
say 
$$\{(U_n\cup X_n:n<\om\}\sq \uu$$  
where each $U_n$ open in $M$ and $X_n\sq X$. 
But since $X$ is Borel concentrated on $Y$
we have that $X\sm\cup\{U_n:n<\om\}$
is countable, so we need only add countably many more
elements of $\uu$ to cover all of $M(X)$.

\bigskip\noindent $M(X)$ is the union of two metrizable subspaces.  
Let
\par $M_1=(X^2\sm \diag)\cup Y^2$ and
\par $M_2= (X\times Y)\cup (Y\times X)\cup (X^2\cap \diag)$.

\par\noindent 
Note that $M_1$ is $N(X^2\sm \diag)$ where 
$N=(X^2\sm \diag)\cup Y^2$ in its separable metric
topology as a subspace of $2^\om\times 2^\om$. By the Lemma
we have that  $X^2\sm \diag$ is relatively $F_\sig$ in $N$ and
so by Proposition \ref{prop} $M_1$ is metrizable. 

To see that $M_2$ is
metrizable use the Bing Metrization Theorem:

\begin{quote}
A topological space is metrizable iff it is regular and
has a $\sig$-discrete basis.
\end{quote}

A family $\bb$ of subsets of $X$ is discrete iff every point of $X$ has
a neighborhood meeting at most one element of $\bb$.
$\sig$-discrete means the countable union of discrete families.

Note that for each $x\in X$ the sets $\{x\}\times Y$ and $Y\times\{x\}$ are
open in $M_2$. Let $\bb$ be a countable open basis for $Y$.
Then
$$\CC=\{U\times\{x\},\;\;\{x\}\times U,\;\; 
\{(x,x)\}\;: \;\; x\in X, U\in\bb\}$$
is an open basis for $M_2$. It is $\sig$-discrete.
The family $\{\{(x,x)\}:x\in X\}$ is
discrete in $M_2$ since $X^2\cap \diag$ is
closed in $M_2$. And for each fixed $U\in\bb$ 
the family $\{U\times\{x\}\;:\;x\in X\}$ is discrete in $M_2$.
(For $(x,x)\in X$ use the neighborhood $\{x\}\times \{x\}$.
For $(y,x)$ with $y\in Y$ and $x\in X$ use the neighborhood $Y\times\{x\}$
and for $(x,y)$ use the neighborhood $\{x\}\times Y$.)  Similarly,
for each $U\in\bb$
the family $\{\{x\}\times U\;:\;x\in X\}$ is discrete in $M_2$.
Since $\bb$ is countable,
$M_2$ has a $\sig$-discrete
basis and is therefor metrizable.

This proves Theorem \ref{one}.
\qed

The next Theorem is an easy generalization of Theorem \ref{one}
using the tower cardinal $\tt$
which is defined as follows.  $\tt$ is the minimum cardinality
of a set $T\sq\infsets$ which is linearly ordered by 
$\sq^*$ but there does not exist $z\in\infsets$ with
$z\sq^* y$ for every $y\in T$.
Martin's axiom implies that $\tt=\cc$.  

\bigskip
\begin{theorem}
Suppose $\tt=\cc$. Then
there exists a nonmetrizable 
regular paracompact space $X$ such that $X^2$ is the union of 
two metrizable subspaces.
\end{theorem}

\proof
The main Lemma changes to:

\begin{lemma}
($\tt=\cc$) There are
 disjoint sets $X,Y\sq 2^\om$ of cardinality $\cc$
such that 
\begin{enumerate}
 
 \item  $X$ is Borel $\cc$-concentrated on
 $Y$, i.e., for every Borel set $B$ in $2^\om$,
 if $Y\sq B$ then $|X\sm B|<\cc$,

 \item $Y^2\sm\diag$ is $F_\sig$ in $X^2\cup Y^2$, and
 
 \item $X^2\sm\diag$ is $F_\sig$ in $X^2\cup Y^2$.

\end{enumerate}
\end{lemma}

The proof is similar.  The space $M=X\cup Y$ is the same.   
Since $X$ is not relatively Borel in $M$ we have by
Proposition \ref{prop} that $M(X)$ is not metrizable.   
But $M(X)$ is regular and paracompact for any $X\sq M$
and metric $M$, see example 5.1.22 Engelking \cite{engel}.
\qed

\bigskip
Remark. The Michael line is the topological
space $M(X)$ where $M$ is  the unit interval, $[0,1]$,
and $X$ the irrationals in $[0,1]$.
Michael Granado in unpublished work has shown that the 
square of the Michael line is not the union 
of two metrizable subspaces.

\begin{question}
(Using just ZFC) Do there exist disjoint sets 
of reals $X$ and $Y$ such that $X$ is not $F_\sig$ in $X\cup Y$ but
$X^2\sm\diag$ is $F_\sig$ in  $X^2\cup Y^2$?
\end{question}

\begin{flushleft}
Arnold W. Miller \\
miller@math.wisc.edu \\
http://www.math.wisc.edu/$\sim$miller\\
University of Wisconsin-Madison \\
Department of Mathematics, Van Vleck Hall \\
480 Lincoln Drive \\
Madison, Wisconsin 53706-1388 \\
\end{flushleft}

\appendix

\newpage

The appendix is not intended for final publication but for  
the electronic version only.

\begin{center}
Appendix 
\end{center}

\begin{quote}
Suppose $M$ is a metric space and $X\sq M$.  Then $M(X)$ is metrizable
iff $X$ is an $F_\sig$ in $M$. (Engelking 5.5.2)
\end{quote}

\proof 
Suppose $X$ is not $F_\sig$ in $M$, then $Y=M\sm X$ is closed 
in $M(X)$ ( since the points of $X$ are isolated, $X$ is open ). 
But $Y$ is not $G_\del$ in $M(X)$.   
To see, this suppose that $Y=\cap_{n\in\om} U_n$ were 
each $U_n$ is open in $M(X)$.  Then there would exists
$V_n$ open in $M$ and $X_n\sq X$ with
$U_n=V_n\cup X_n$. But then $Y= \cap_{n\in\om} V_n$ which contradicts
$Y$ is not $G_\del$ in $M$.

For the converse,
suppose $X$ is $F_\sig$ in $M$ and write it as the 
union of closed sets $X=\cup_{n<\om}C_n$.  $M(X)$ is regular
so it is enough by the Bing Metrization Theorem to check that it 
has a $\sig$-discrete base.  Let $\bb$ 
be a $\sig$-discrete base for $M$. We claim that 
$$\bb \cup\{\{x\}:x\in X\}$$
which is a basis for $M(X)$ is $\sig$-discrete in $M(X)$.  
$\bb$ is $\sig$-discrete in $M$ so it is also $\sig$-discrete in $M(X)$.
$$\{\{x\}:x\in X\}=\cup_{n<\om}\CC_n\mbox{ where } \CC_n=\{\{x\}:x\in C_n\}$$
shows that it is $\sig$-discrete, since for any $n$ if
$x\notin C_n$ then $M\setminus C_n$ is a neighborhood of $x$ missing
all elements of $\CC_n$.

\bigskip
\begin{quote}
$M(X)$ is regular paracompact, whenever $M$ is metric. (Engelking 5.1.22)
\end{quote}
\proof
Regular: Given $p\in M$ if $p\in X$ then it is has the clopen neighborhood 
$\{p\}$, if $p\notin X$, then 
the neighborhoods of $p$ in $M$ are also a neighborhood basis in $M(X)$.

Paracompact: Let $\uu$ be an open cover of basic open sets in $M(X)$.  We 
may assume it has the form: 
 $$\uu=\vv\cup \{\{x\}:x\in Z\}$$
where $\vv$ is a family of basic open sets in $M$ and  $Z=X\sm \cup\vv$. Since
metric spaces are hereditarily paracompact, there exists a locally finite
refinement $\ww$ of $\vv$ with $\cup\vv=\cup\ww$. 
But then $\ww\cup \{\{x\}:x\in Z\}$ is a 
locally finite refinement of $\uu$. 

\end{document}